\documentclass[reqno]{amsart}

\usepackage{latexsym}
\usepackage{amsmath}
\usepackage{amssymb}
\numberwithin{equation}{section}

\begin{document}


{\theoremstyle{plain}%
  \newtheorem{theorem}{Theorem}[section]  
  \newtheorem{corollary}[theorem]{Corollary}
  \newtheorem{proposition}[theorem]{Proposition}
  \newtheorem{lemma}[theorem]{Lemma}
  \newtheorem{question}[theorem]{Question}
  \newtheorem{conjecture}[theorem]{Conjecture}
}

{\theoremstyle{remark}
  \newtheorem{fact}{Fact}
  \newtheorem{construction}[theorem]{Construction}
}

{\theoremstyle{definition}
  \newtheorem{definition}[theorem]{Definition}
  \newtheorem{remark}[theorem]{Remark}
  \newtheorem{example}[theorem]{Example}
}

\newcommand{\m}[1]{\marginpar{\addtolength{\baselineskip}{-3pt}{\footnotesize \it #1}}}
\newcommand{\reg}{\operatorname{reg}}
\newcommand{\depth}{\operatorname{depth}}
\newcommand{\pdim}{\operatorname{pd}}
\newcommand{\lcm}{\operatorname{lcm}}
\newcommand{\lk}{\operatorname{conn}}
\newcommand{\rlk}{\overline{\operatorname{conn}}}
\newcommand{\M}{\mathcal{M}}
\newcommand{\I}{\mathcal{I}}
\newcommand{\F}{\mathcal{F}}
\newcommand{\A}{\mathcal{A}}
\newcommand{\K}{\mathcal{K}}
\newcommand{\G}{\mathcal{G}}
\newcommand{\N}{\mathbb{N}}
\newcommand{\pr}{\mathbb{P}}


\title{Splittable ideals and the resolutions of monomial ideals}
\thanks{Version: \today}
\author{Huy T\`ai H\`a}
\address{Tulane University \\ Department of Mathematics \\
6823 St. Charles Ave. \\ New Orleans, LA 70118, USA}
\email{tai@math.tulane.edu}
\urladdr{http://www.math.tulane.edu/$\sim$tai/}

\author{Adam Van Tuyl}
\address{Department of Mathematical Sciences \\
Lakehead University \\
Thunder Bay, ON P7B 5E1, Canada}
\email{avantuyl@sleet.lakeheadu.ca}
\urladdr{http://flash.lakeheadu.ca/$\sim$avantuyl/}


\begin{abstract}
We provide a new combinatorial approach to
study the minimal free resolutions of edge ideals,
that is, quadratic square-free monomial ideals.
With this method
we can recover most of the known results on resolutions of
edge ideals with fuller generality, and at the same time, obtain new results.
Past investigations on the resolutions of edge ideals usually reduced the problem
to computing the dimensions of reduced homology or Koszul homology groups.
Our approach circumvents the highly nontrivial problem of computing the 
dimensions of these groups and turns
the problem into combinatorial questions about the associated simple graph.
We also show that our technique extends
successfully to the study of graded Betti numbers of arbitrary 
square-free monomial ideals viewed as facet
ideals of simplicial complexes.
\end{abstract}

\maketitle


\section{Introduction}

Let $R = k[x_1,\ldots,x_n]$ be a polynomial ring over an arbitrary field $k$.
If $I$ is a homogeneous ideal
of $R$, then associated to $I$ is a {\it minimal
graded free resolution}
\[
0 \rightarrow \bigoplus_j R(-j)^{\beta_{l,j}(I)}
\rightarrow \bigoplus_j R(-j)^{\beta_{l-1,j}(I)}
\rightarrow \cdots
\rightarrow  \bigoplus_j R(-j)^{\beta_{0,j}(I)}
\rightarrow I \rightarrow 0
\]
where the maps are exact, $l \leq n$, and $R(-j)$ is the
$R$-module shifted by $j$.
The number $\beta_{i,j}(I)$, the $ij$th
{\it graded Betti number} of $I$, is an invariant of $I$ 
equal to the number of minimal generators of degree $j$ in the $i$th syzygy module.

In this paper we shall study the graded Betti numbers of {\it monomial ideals}. The book of Miller and Sturmfels \cite{MS} contains a comprehensive introduction and 
list of references on this topic. We shall concentrate on ideals
which are generated by square-free quadratic monomials
so that we may exploit the natural bijection
\[\left\{
\begin{array}{c}
\mbox{square-free quadratic monomial} \\
\mbox{ideals $I \subseteq R = k[x_1,\ldots,x_n]$}
\end{array}
\right\} \leftrightarrow
\left\{
\mbox{simple graphs $G$ on $n$ vertices}
\right\}.
\]
By a simple graph we mean an undirected graph with no loops or multiple edges,
but not necessarily connected.  The bijection is defined by mapping
the graph  $G$ with edge set $E_G$ and vertices $V_G = \{x_1,\ldots,x_n\}$
to the square-free monomial ideal
\begin{equation}\label{IGdefn}
\I(G) = ( \{x_ix_j ~|~ \{x_i,x_j\} \in E_G\}) \subseteq k[x_1, \dots, x_n].
\end{equation}
(The graph of $n$ isolated vertices is mapped to $I = (0)$ which we shall
also consider as a square-free quadratic monomial ideal.) Note that
(\ref{IGdefn}) implies that $\beta_{0,j}(\I(G)) = |E_G|$ if $j = 2$,
and  $0$ if $j \neq 2$.  Thus, we shall be interested in formulas for $\beta_{i,j}(\I(G))$
with $i \geq 1$.
The ideal $\I(G)$ is commonly called the {\it edge ideal} of $G$.

Edge ideals, which were first introduced
by Villarreal \cite{V}, are the focus of an ongoing program in algebraic combinatorics.
Many authors have been interested in establishing a dictionary between
the algebraic invariants of $\I(G)$ and the combinatorial data associated to the
graph $G$.
The references \cite{eghp,EV,FV,FH, Fr,HHZ,HHZ1,J,JK,K,RVT,SVV,V,V1,V2,Z} form a partial
list of references on this topic.

Because the edge ideal is a square-free monomial ideal, Hochster's formula
\cite{Ho} and its variant, the
Eagon-Reiner formula \cite{ER}, provide us with tools to study the 
numbers $\beta_{i,j}(\I(G))$.
More precisely, the ideal $\I(G)$ can be associated with a
simplicial complex $\Delta(G)$ via the Stanley-Reisner correspondence.
The numbers $\beta_{i,j}(\I(G))$ are then
related to the dimensions of the reduced homology groups of subcomplexes
of $\Delta(G)$, or as in the case of Eagon and Reiner's formula, the Alexander dual
of $\Delta(G)$.  An examination of the papers \cite{eghp,Fr,J,JK,K,RVT,Z}
reveals that these formulas provide the basis for most of the known results on the numbers
$\beta_{i,j}(\I(G))$.
The exception to this observation is \cite{Z}
which uses Koszul homology.

In this paper we use the notion of a {\it splittable monomial ideal}, as
first defined by Eliahou and Kervaire \cite{EK}, to
introduce a new technique to the study the numbers $\beta_{i,j}(\I(G))$.
Our approach, which has
the advantage that we can avoid the highly nontrivial problem of computing
the dimensions of reduced homology or Koszul homology groups, allows us to 
recover many of the known
results with fuller generality, and at the same time, provides new results.  The use of
splittable monomial ideals also provides a unified combinatorial perspective for most of
the known results.

A monomial ideal $I$ is {\it splittable} if there exists two monomial ideals $J$ and $K$
such that $I = J + K$, and furthermore, the generators of $J \cap K$ satisfy certain
technical conditions (see
Definition \ref{defn: split}).
Splittable ideals allow us to
relate $\beta_{i,j}(I)$ to the graded Betti numbers of the ``smaller'' ideals 
$J, K$ and $J \cap K$
(see Theorem \ref{prop: ekf}).
Given an edge ideal $\I(G)$, our goal is to find a splitting $\I(G) = J + K$
so that $J, K$, and $J \cap K$ are related to edge ideals of subgraphs
of $G$, and therefore produce a recursive like formula for the graded Betti numbers of $\I(G)$.

Sections 3 and 4 of the paper are devoted to  two natural candidates for a splitting of 
$\I(G)$.  If $e = \{u,v\}$ is any edge of $G$, then it is clear that
\begin{equation}\label{split1}
\I(G) = (uv) + \I(G\backslash e)
\end{equation}
where $G \backslash e$ is the subgraph of $G$ with the edge $e$ removed.  Similarly,
if $v$ is any vertex of $G$, and if $N(v) = \{v_1,\ldots,v_d\}$ denotes the
distinct neighbors of $v$, then
\begin{equation}\label{split2}
\I(G) = (vv_1,vv_2,\ldots,vv_d) + \I(G \backslash\{v\})
\end{equation}
where $G \backslash \{v\}$ is the subgraph of $G$ with vertex $v$ and edges 
incident to $v$ removed. Observe
that $(vv_1,\ldots,vv_d)$ is the edge ideal of the complete bipartite graph $\K_{1,d}$.

In general (\ref{split1}) and (\ref{split2}) will not be splittings of $\I(G)$. We therefore
call $e$ a {\it splitting edge} of $G$ if  (\ref{split1}) is
a splitting, and similarly, we say $v$ is a {\it splitting vertex} if
(\ref{split2}) is a splitting.
Theorems  \ref{classify splitting edges} and \ref{prop: edgeidealssplit} then
characterize
which edges and vertices of $G$ can have this property.  An edge $e = \{u,v\}$
is a splitting edge if  the set of  neighbors of $u$ (or $v$) is a 
subset of  $N(v) \cup \{v\}$ (or
$N(u) \cup \{u\}$).  Splitting vertices are more ubiquitous; a vertex $v$ is a splitting
vertex provided $v$ is not an isolated vertex or the vertex of degree $d$ of $\K_{1,d}$.

We adopt the convention that for any ideal $I$,  $\beta_{-1,j}(I) = 1$ if $j=0$, and  
$\beta_{-1,j}(I) = 0$ otherwise. Our first main result is the following formulas for $\beta_{i,j}(\I(G))$:

\begin{theorem} \label{intro main theorem}
Let $G$ be a simple graph with edge ideal $\I(G)$.

\noindent
$(i)$ {\rm (Theorem \ref{recursive formula 3})}
Suppose $e = \{u,v\}$ is a splitting edge of $G$. Set $H = G \backslash (N(u) \cup N(v))$.
If $n = |N(u) \cup N(v)| - 2$, then for all $i \geq 1$
\[\beta_{i,j}(\I(G)) = \beta_{i,j}(\I(G \backslash e)) +
\sum_{l=0}^i \binom{n}{l} \beta_{i-l-1,j-2-l}(\I(H)).\]
\noindent
$(ii)$ {\rm (Theorem \ref{recursive formula 2})}
Let $v$ be a splitting vertex  of $G$ with $N(v) = \{v_1,\ldots,v_d\}$.
Set
\[G_i :=  G \backslash (N(v) \cup N(v_i)) ~~\mbox{for $i = 1,\ldots,d,$}\] and let
$G_{(v)}$ be the subgraph of $G$ consisting of all edges incident to $v_1,\ldots,v_d$
except those that are also incident to $v$.
Then for all $i,j \geq 0$
\[\beta_{i,j}(\I(G)) = \beta_{i,j}(\I(\K_{1,d})) + \beta_{i,j}(\I(G \backslash \{v\})) +
\beta_{i-1,j}(L)\]
where  $L = v\I(G_{(v)}) + vv_1\I(G_1) + \cdots +vv_d\I(G_d)$.
\end{theorem}

Our formula in Theorem \ref{intro main theorem} $(i)$ unifies all known results about
$\beta_{i,j}(\I(G))$ when $G$ is a forest. Since a leaf is a splitting edge,
we recover the recursive formula for the graded Betti numbers of forests (cf.
Corollary \ref{treeformula}) as first given by Jacques and Katzman \cite{J,JK}.  
At the same time our recursive formula
is more general since it applies to all leaves, and
not only the special leaf required in the argument of \cite{J,JK}.
Theorem \ref{intro main theorem} $(i)$
also allows us to give a combinatorial proof (cf. Corollary \ref{tree reg}) 
of Zheng's formula \cite{Z}
for the regularity of $\I(G)$ in terms of the number of disconnected edges in  $G$
when $G$ is a forest.
The above formula fails to be recursive in 
general because the subgraphs may not contain splitting edges, 
thus preventing
us from reiterating the process. However, it is still general enough to provide new 
results on the projective dimension and regularity of edge ideals 
(cf. Corollary \ref{reg pdim 2}).

The formula in Theorem \ref{intro main theorem} $(ii)$ is not recursive because it
involves computing $\beta_{i,j}(L)$ where $L$ is not an edge ideal. Yet, this 
formula proves to be
very effective in studying the {\it linear strand} of the minimal free 
resolution of edge ideals.
We can give new combinatorial proofs for many of the results on the linear 
strand of edge ideals;
for example, we can recover (cf. Corollary \ref{np property}) a result of 
Eisenbud, et al. \cite{eghp}
on the $N_{2,p}$
property of edge ideals, a notion closely tied to the $N_p$ {\it property} introduced
by Green \cite{Gr}. Precisely,
Theorem \ref{intro main theorem} $(ii)$ enables us to give
a new proof to the fact that $\I(G)$ has property $N_{2,p}$ for $p > 1$
if and only if every minimal cycle of $G^c$, the complementary graph of $G$, has 
length at least
$p+3$.  Theorem \ref{intro main theorem} $(ii)$ also allows us to recover the
formulas for the Betti numbers in the linear strand first shown in \cite{RVT} 
(cf. Corollary \ref{linear strand}).
At the same time, we can provide new results on the projective dimension and regularity of edge
ideals (cf. Corollary \ref{reg pdim 3}).

In the final section we extend the scope of this paper by considering
the graded Betti numbers of facet ideals.
The facet ideal was
introduced by Faridi \cite{faridi:2002,faridi:2004} to generalize an edge ideal. 
Since any square-free
monomial ideal can be realized as the facet ideal of a simplicial complex, 
our method thus works for a
large class of monomial ideals.
Let $\Delta$ be a simplicial complex on the vertex set 
$V_\Delta = \{ x_1, \dots, x_n\}$. The {\it facet ideal} 
$\I(\Delta)$ of $\Delta$ is defined to be
$$\I(\Delta) = ( \{ \prod_{x \in F} x ~|~ F \mbox{ is a facet of } \Delta \} ) 
\subseteq k[x_1, \dots, x_n].$$
If $F$ is a facet of $\Delta$, then we say that $F$ is a {\it splitting facet} if
$\I(\Delta) = (\prod_{x \in F} x) + \I(\Delta')$,
where $\Delta' = \Delta \backslash F$ is the subcomplex of $\Delta$ with the facet $F$ removed,
is a splitting of $\I(\Delta)$. Our next main result is a higher dimension analogue of 
Theorem \ref{intro main theorem} $(i)$ relating the graded Betti number of $\I(\Delta)$ 
to those of facet ideals of subcomplexes of $\Delta$; 
see Definition \ref{3conn} for unexplained terminology.

\begin{theorem}[Theorem \ref{simplicial:betti}] \label{intro main theorem 2}
Let $F$ be a splitting facet of a simplicial complex $\Delta$. Let 
$\Delta' = \Delta \backslash F$ and $\Omega = \Delta \backslash \lk_\Delta(F).$ 
Then for all $i \geq 1$ and $j \ge 0$,
$$\beta_{i,j}(\I(\Delta)) = \beta_{i,j}(\I(\Delta')) + \sum_{l_1=0}^i
\sum_{l_2=0}^{j-|F|} 
\beta_{l_1-1,l_2}(\I(\rlk_\Delta(F))) \beta_{i-l_1-1,j-|F|-l_2}(\I(\Omega)).$$
\end{theorem}

Similar to Theorem \ref{intro main theorem}, our formula in 
Theorem \ref{intro main theorem 2} is recursive for simplicial forests 
(cf. Theorems \ref{3split} and \ref{simplicial:recursive}).  Consequently, there exists
a large class of square-free monomial ideals that can be examined via a recursive formula.  
The recursive formula of Jacques and Katzman for forests \cite{JK} becomes a special case 
of this result.  Moreover, formulas for the graded Betti numbers in the linear strand 
of facet ideals of simplicial forests are recovered, generalizing results of Zheng \cite{Z}.

\noindent
{\bf Acknowledgments.}
The computer algebra package
{\tt CoCoA} \cite{Co} was used extensively to generate examples.
The authors would like to thank  Hema Srinivasan for her comments.
The second author was partially supported by NSERC.


\section{Preliminaries}

For completeness we gather together the needed results and definitions on simple graphs, simplicial complexes,
resolutions, and splittable ideals.  Readers familiar with this material may wish to continue
directly to the next section.

\subsection{Graph terminology, simplicial complexes, edge and facet ideals}
In this paper $G$ will denote a finite simple graph (undirected, no
loops or multiple edges, but not necessarily connected). We denote by $V_G$ and $E_G$ 
the set of vertices
 and edges, respectively, of $G$.

If $V_G = \{x_1, \dots, x_n\}$, then we associate to $G$ a polynomial ring 
$R = k[x_1, \dots, x_n]$
 (here, by abuse of notation, we use the $x_i$s to denote both the vertices in 
$V_G$ and the variables
 in the polynomial ring). For simplicity we write $uv \in E_G$ instead of $\{u,v\} \in E_G$.
Also, by abuse of notation, we use $uv$ to denote both the edge $uv$ and the monomial $uv$ in
 the edge ideal. In particular,
$\I(G) = ( \{ uv ~|~ uv \in E_G \} ) \subseteq R.$

A vertex $y$ is a {\it neighbor} of $x$ if $xy \in E_G$.  Set $N(x) :=
\{y \in V_G ~|~ xy \in E_G\}$, the set of all neighbors
of $x$ in $G$. The {\it degree} of a vertex
$x \in V_G$, denoted by $\deg_G x$, is the number of edges incident to $x$. When there is no
confusion, we shall omit $G$ and write $\deg x$. Observe that $\deg x = |N(x)|$ since $G$ is simple.

If $e \in E_G$, we shall write $G\backslash e$ for the subgraph of $G$ with
the edge $e$ deleted.  If $S = \{x_{i_1},\ldots,x_{i_s}\} \subseteq V_G$,
we shall write $G \backslash S$ for the subgraph of $G$ with the vertices
of $S$ (and their incident edges) deleted.  We further write $G_S$ to denote the
{\it induced subgraph} of $G$ on $S$ (i.e., the subgraph of $G$ whose vertex set is 
$S$ and whose edges
are edges of $G$ connecting vertices in $S$). We say that $C = (x_1x_2 \dots x_lx_1)$ 
is a {\it cycle}
 of $G$ if $x_ix_{i+1} \in E_G$ for $i = 1, \dots, l$ (where $x_{l+1} = x_1$).
The {\it complete} graph $\K_n$ of size $n$ is the graph whose vertex set $V$ has $n$ 
vertices and whose
edges are $\{ uv ~|~ u \neq v \in V \}$. A complete graph $\K_n$ which is a subgraph of 
$G$ is
called a {\it $n$-clique} of $G$. The {\it complete bipartite} graph $\K_{m,n}$ is the 
graph whose vertex
set can be divided into two disjoint subsets $A$ and $B$ such that $|A| = m, |B| = n$,
and the edges of the graph are $\{ uv ~|~ u \in A, v \in B \}$.

A {\it simplicial complex} $\Delta$ over a vertex set $V_\Delta = \{ x_1, \dots, x_n\}$
is a collection of subsets of $V_\Delta$, with the property that $\{x_i\} \in \Delta$ for
all $i$, and if $F \in \Delta$ then all subsets of $F$ are also in $\Delta$.
Elements of $\Delta$ are called {\it faces}.
The {\it dimension} of a face $F$, denoted by $\dim F$, is defined to
 be $|F|-1$, where $|F|$ denotes the cardinality of $F$. The {\it dimension}
 of $\Delta$, denoted by $\dim \Delta$, is defined to be the maximal dimension of a
 face in $\Delta$. The maximal faces of $\Delta$ under inclusion are called {\it facets}. 
If all facets of $\Delta$ have the same dimension $d$, then $\Delta$ is
said to be {\it pure $d$-dimensional}.

We usually denote the simplicial complex $\Delta$ with facets $F_1, \dots, F_q$ by
$\Delta = \langle F_1, \dots, F_q \rangle;$
here, the set $\F(\Delta) = \{ F_1, \dots, F_q \}$ is often referred to as the
{\it facet set} of $\Delta$. If $F$ is a facet of $\Delta$, say $F = F_q$, then
 we denote by $\Delta \backslash F$ the simplicial complex obtained by removing
 $F$ from the facet set of $\Delta$, i.e., $\Delta \backslash F = \langle F_1, \dots, F_{q-1} \rangle.$
Throughout the paper, by a {\it subcomplex} of a simplicial complex $\Delta$, we
shall mean a simplicial complex whose facet set is a subset of the facet set of
$\Delta$. If $\Delta'$ is a subcomplex of $\Delta$, then we denote by
$\Delta \backslash \Delta'$ the simplicial complex obtained from $\Delta$ by
 removing from its facet set all facets of $\Delta'$.

We say that two facets $F$ and $G$ of $\Delta$ are {\it connected} if there exists a
{\it chain} of facets of $\Delta$, $F = F_0, F_1, \dots, F_m = G$, such that
$F_i \cap F_{i+1} \not= \emptyset$ for any $i = 0, \dots, m-1$. The simplicial
complex $\Delta$ is said to be {\it connected} if any two of its facets are connected.

To a simplicial complex $\Delta$ over the vertex set $V_\Delta = \{ x_1, \dots, x_n \}$
we associate an ideal $\I(\Delta)$ in the
polynomial ring $R = k[x_1, \dots, x_n]$.
We write $F$ to denote both a facet of $\Delta$ and the monomial $\prod_{x \in F} x$.
In particular,
$\I(\Delta) = ( \{ F ~|~ F \in \F(\Delta) \} ) \subseteq R.$

A facet $F$ of $\Delta$ is a {\it leaf}
if either $F$ is the only facet of $\Delta$, or there exists a facet
$G$ in $\Delta$, $G \not= F$, such that $F \cap F' \subseteq F \cap G$ for every
facet $F' \in \Delta, F' \not= F$. The simplicial complex $\Delta$ is
called a {\it tree} if $\Delta$ is connected and every nonempty connected
subcomplex of $\Delta$ (including $\Delta$ itself) has a leaf. We call $\Delta$ a
{\it forest} if every connected component of $\Delta$ is a tree.

\subsection{Resolutions, Betti numbers, and splittable ideals}
Let $\G(I)$ denote the minimal set  of generators of a monomial ideal $I$; this set is uniquely
determined (cf. Lemma 1.2 of \cite{MS}). The following definition and result play an essential role throughout the paper.

\begin{definition}[see \cite{EK}]\label{defn: split}
A monomial ideal $I$ is {\it splittable} if $I$ is the sum
of two nonzero monomial ideals $J$ and $K$, that is, $I = J+K$, such
that
\begin{enumerate}
\item $\G(I)$ is the disjoint union of $\G(J)$ and $\G(K)$.
\item there is a {\it splitting function}
\begin{eqnarray*}
\G(J\cap K) &\rightarrow &\G(J) \times \G(K) \\
w & \mapsto & (\phi(w),\psi(w))
\end{eqnarray*}
satisfying
\begin{enumerate}
\item for all $w \in \G(J \cap K), ~~ w = \lcm(\phi(w),\psi(w))$.
\item for every subset $S \subset \G(J \cap K)$, both
$\lcm(\phi(S))$ and $\lcm(\psi(S))$
strictly divide $\lcm(S)$.
\end{enumerate}
\end{enumerate}
If $J$ and $K$ satisfy the above properties, then
we say $I = J + K$ is a {\it splitting} of $I$.
\end{definition}

\begin{theorem}[Eliahou-Kervaire \cite{EK} Fatabbi \cite{Fa}]
\label{prop: ekf}
Suppose $I$ is a
splittable monomial ideal with splitting $I = J+K$.  Then
for all $i, j \geq 0$,
\[\beta_{i,j}(I) = \beta_{i,j}(J) + \beta_{i,j}(K) +
\beta_{i-1,j}(J\cap K).\]
\end{theorem}

Recall that for an ideal $I$ generated by elements of degree at least $d$, the Betti numbers
$\beta_{i,i+d}(I)$ form the so-called {\it linear strand} of $I$ (see \cite{eghp,hi}). An  ideal $I$ generated
by elements of degree $d$ is said to have a {\it linear resolution}
if the only nonzero graded Betti numbers are those in the linear strand. 
Of particular interest are also the following invariants which
measure the ``size'' of the  minimal graded free resolution of $I$.
The {\it regularity} of $I$,  denoted $\reg(I)$, is defined
by
\[\reg(I) := \max\{j-i ~|~ \beta_{i,j}(I) \neq 0\}.\]
The {\it projective dimension} of $I$, denoted $\pdim(I)$, is defined to be
\[\pdim(I):= \max\{i ~|~ \beta_{i,j}(I) \neq 0 \}.\]

When $I$ is a splittable ideal, Theorem \ref{prop: ekf} implies the following result:

\begin{theorem} \label{reg pdim}
If $I$ is a splittable monomial ideal with splitting $I = J+K$, then
\begin{enumerate}
\item[$(i)$] $\reg(I) = \max\{\reg(J),\reg(K),\reg(J\cap K) - 1\}$.
\item[$(ii)$] $\pdim(I) = \max\{\pdim(J),\pdim(K),\pdim(J\cap K) + 1\}$.
\end{enumerate}
\end{theorem}

The following results shall be required throughout the paper.  The
lemma is well known.   See, for example, Lemma 2.1 and Corollary 2.2 of \cite{JK}.

\begin{lemma}\label{betti-tensor}
Let $R = k[x_1,\ldots,x_n]$ and $S = k[y_1,\ldots,y_m]$, and
let $I \subseteq R$ and $J \subseteq S$ be homogeneous ideals.  Then
\[\beta_{i,j}(R/I \otimes S/J) = \sum_{l_1 = 0}^i \sum_{l_2=0}^j
\beta_{l_1,l_2}(R/I)\beta_{i-l_1,j-l_2}(S/J).\]
\end{lemma}

\begin{remark} \label{rmk-tensor}
If $I, J \subseteq R = k[x_1, \dots, x_n]$ are square-free monomial ideals such 
that none of the $x_i$s appearing in the minimal generators of $I$ appear in 
the minimal generators of $J$,  then $R/I \otimes R/J = R/(I+J)$. 
Lemma \ref{betti-tensor} thus implies
$$\beta_{i,j}(R/(I+J)) = \sum_{l_1 = 0}^i \sum_{l_2=0}^j \beta_{l_1, l_2}(R/I) \beta_{i-l_1, j-l_2}(R/J).$$
\end{remark}

\begin{remark}\label{betti-non-con}
If $G$ is a simple graph with two connected components, i.e.,
$G = G_1 \cup G_2$, with $V _G = V_{G_1} \cup  V_{G_2}$ and $V_{G_1} \cap V_{G_2} = \emptyset$,
then Remark \ref{rmk-tensor} implies that to calculate $\beta_{i,j}(\I(G))$, it is enough 
to calculate the graded Betti numbers of the edge ideals  $\I(G_1)$ and $\I(G_2)$.  
More generally, if $G$ has $n \geq 2$ components, by repeated applying Remark 
\ref{rmk-tensor}, to calculate  $\beta_{i,j}(\I(G))$ it suffices to
 calculate the Betti numbers of the edge ideals associated to each connected component of $G$.
\end{remark}

\begin{theorem}\label{resk1d}
Suppose that $G = \K_{1,d}$.  Then for $i \geq 0$
\[\beta_{i,j}(\I(G)) =
\left\{
\begin{array}{ll}
\binom{d}{i+1} & \mbox{if $j=i+2$} \\
0 & \mbox{otherwise.}
\end{array}
\right.\]
\end{theorem}
\begin{proof}
Since $G = \K_{1,d}$, it follows that
$\I(G) = (vv_1,\ldots,vv_d) \subseteq R = k[v,v_1,\ldots,v_d].$
The conclusion now follows from the fact that $v_1,\ldots,v_d$ is
a regular sequence on $R$, and that $\beta_{i,j}(\I(G)) = \beta_{i,j-1}((v_1,\ldots,v_d))$.
\end{proof}

\section{Splitting Edges}

Let $G$ be a simple graph with edge ideal $\I(G)$ and  $e = uv \in E_G$.
If we set $J = (uv) ~\mbox{and}~  K= \I(G\backslash e),$
then $\I(G) = J + K$.  In general this
may not be a splitting of $\I(G)$.  The goal of this section is to determine
when $J$ and $K$ give a splitting of $\I(G)$, and furthermore, how
this splitting can be used to ascertain information about the numbers $\beta_{i,j}(\I(G))$.

We begin by assigning a name to an
edge for which there is a splitting.

\begin{definition}An edge $e = uv$ is a
{\it splitting edge} if $\I(G) = (uv) + \I(G\backslash e)$ is a splitting.
\end{definition}

\begin{lemma} \label{Ideal J K}
Let $J = (uv)$ and $K = \I(G\backslash e)$ with $e = uv \in E_G$.   
If $N(u) \backslash \{v\}
= \{u_1,\ldots,u_n\}$, $N(v)
\backslash \{u\} = \{v_1,\ldots,v_m\}$, and $H = G \backslash (N(u) \cup N(v))$, then
\[ J \cap K = uv((u_1,\ldots,u_n,v_1,\ldots,v_m) + \I(H)).\]
\end{lemma}

\begin{proof}
Because $J = (uv)$ and $K = \I(G\backslash e)$ are both monomial ideals,
\[ J \cap K = ( \{\operatorname{lcm}(uv,m) ~|~ m \in \G(K)\} ).\]
Each $m \in  \G(K)$ corresponds to an edge of $G\backslash e$. There
are three cases for this edge:  (1) it is
incident to either $u$ or $v$, (2) it is not incident to $u$ or $v$, but is
incident to a neighbor of either $u$ or $v$, or (3) it is not incident
to any vertex in $N(u) \cup N(v)$.

If $m$ is in  cases (1) and (2),
then $\lcm(uv,m)$ is in $uv(u_1, \dots, u_n, v_1, \dots, v_m)$.  If $m$ is in case (3),
$\lcm(uv,m)$ belongs to $uv\I(H)$. The statement follows.
\end{proof}

We in fact obtain the following description for $\G(J \cap K)$.
\begin{corollary} \label{UV gens}
Let $e = uv \in E_G, J = (uv)$ and $K = \I(G\backslash e)$.
If $A = N(u) \backslash \{v\}$ and $B = N(v)\backslash \{u\}$, then
\begin{eqnarray*}
 \G(J \cap K)  &=&  \{ uvu_i ~|~ u_i \in A\backslash B\} \cup 
\{ uvv_i ~|~ v_i \in B\backslash A\} \cup  \{ uvz_i ~|~ z_i \in A \cap B \}
\cup \\
&& \{ uvm ~|~ m \in \I(H)\}.
\end{eqnarray*}
\end{corollary}
\normalsize

The above description of $\G(J \cap K)$ will enable us to identify splitting edges.

\begin{theorem}\label{classify splitting edges}
An edge $e = uv$ is a splitting edge of $G$
if and only if $N(u) \subseteq (N(v) \cup \{v\})$ or
$N(v) \subseteq (N(u) \cup \{u\})$.
\end{theorem}

\begin{proof}
$(\Leftarrow)$.  Without loss of generality, we shall assume
that $N(u) \subseteq (N(v) \cup \{v\})$.  This condition and
Corollary \ref{UV gens} then imply that
\[
\G(J\cap K) = \{ uvv_i ~|~ v_i \in N(v) \backslash \{u\}\} \cup
\{ uvm ~|~ m \in \I(H)\}.\]
To show that $e = uv$ is splitting edge,
it suffices to verify that  the function
$\G(J \cap K) \rightarrow \G(J) \times \G(K)$
defined by
\[w \mapsto (\phi(w),\psi(w)) =
\left\{
\begin{array}{ll}
(uv,vv_i) & \mbox{if $w = uvv_i$} \\
(uv,m) & \mbox{if $w = uvm$}\\
\end{array}
\right.\]
satisfies conditions (a) and (b)
of Definition \ref{defn: split}.  Indeed, condition $(a)$ is immediate.
So, suppose $S \subseteq \G(J \cap K)$.  Our description of
$\G(J \cap K)$ implies all elements of $S$ are divisible by $uv$. Moreover, $\lcm(S)$ must have degree at least three. Thus, $\lcm(\phi(S)) = uv$ strictly divides $\lcm(S)$.
Furthermore, $u$ does not divide $\lcm(\psi(S))$
implying that  $\lcm(\psi(S))$ strictly divides $\lcm(S)$.
Condition $(b)$ now follows.

($\Rightarrow$)
We prove the contrapositive.  Suppose that $e = uv$ is an edge
such that $N(u) \not\subseteq (N(v) \cup \{v\})$ and
$N(v) \not\subseteq (N(u) \cup \{u\})$.  Hence, there exists vertices
$x$ and $y$ such that $ux,vy\in E_G$, but
$uy,vx \not\in E_G$.

We now show that no splitting function can exist.
If there was a splitting function
$\G(J \cap K) \rightarrow \G(J) \times \G(K)$
our splitting function must have the
form
\[w \mapsto (\phi(w),\psi(w)) = (uv,\psi(w))\]
since $\G(J) = \{uv\}$.
By Corollary \ref{UV gens} it follows that $uvx,uvy \in \G(J \cap K)$.
Condition $(a)$ of Definition \ref{defn: split},
would imply $uvx = \lcm(\phi(uvx),\psi(uvx)) =
\lcm(uv,\psi(uvx))$.  Thus $\psi(uvx) = x, vx, ux$ or $uvx$.
But since  $\psi(uvx) \in \G(K)$ and $vx \not\in E_G$,
this forces $\psi(uvx) = ux$.  By a similar argument,
$\psi(uvy) = vy$.

The subset $S = \{uvx,uvy \} \subseteq \G(J\cap K)$
now fails to satisfy Definition \ref{defn: split} (b)
since $\lcm(S) = \lcm(\psi(S)) = uvxy$.
Thus $e = uv$ is not a splitting edge.
\end{proof}

The following identity for the numbers of
$\beta_{i-1,j}(J \cap K)$ can now be derived.

\begin{lemma} \label{betti J cap K}
Let $e = uv$ be a splitting edge of $G$ with
$N(u) \subseteq (N(v) \cup \{v\})$.  If $N(v)\backslash\{v\} = \{v_1,\ldots,v_n\}$,
$J = (uv)$, and $K = \I(G\backslash e)$, then for $i \geq 1$ and all $j \geq 0$
\[\beta_{i-1,j}(J \cap K) =
\sum_{l=0}^i \binom{n}{l} \beta_{i-l-1,j-2-l}(\I(H))\]
where $\I(H)$ is the edge ideal of $H = G\backslash \{u,v,v_1,\ldots,v_n\}$.
\end{lemma}

\begin{proof}
When $e = uv$ is a splitting edge, the conclusion
of Lemma \ref{Ideal J K} becomes
\[J \cap K = uv((v_1,\ldots,v_n) + \I(H)).\]
where $H = G\backslash \{u,v,v_1,\ldots,v_n\}$.
Set $L = (v_1,\ldots,v_n) + \I(H)$.
Since no generator of $L$
is divisible by either $u$ or $v$, we have that $uv$ is a nonzero
divisor on $R/L$.  As a consequence
\[\beta_{i-1,j}(uvL) = \beta_{i-1,j-2}(L) = \beta_{i,j-2}(R/L).\]
 Observe that none of the generators
of $\I(H)$ are divisible by $v_i$ for $i = 1,\ldots,n$. 
Now apply Remark \ref{rmk-tensor} to compute
$\beta_{i,j-2}(R/L)$  and use the
fact that the graded Betti numbers of $R/(v_1,\ldots,v_n)$ are given
by the Koszul resolution.
\end{proof}

We now state and prove the main theorem of this section.

\begin{theorem}\label{recursive formula 3}
Let $e = uv$ be a splitting edge of $G$, and set $H = G \backslash (N(u) \cup N(v))$.  If $n = |N(u) \cup N(v)| - 2$,
then for all $i \geq 1$ and all $j \geq 0$
\[\beta_{i,j}(\I(G)) =  \beta_{i,j}(\I(G \backslash e)) +
\sum_{l=0}^i \binom{n}{l} \beta_{i-l-1,j-2-l}(\I(H)).\]
\end{theorem}

\begin{proof}
Because $e = uv$ is a splitting edge, we can assume without loss
of generality that $N(u) \subseteq (N(v) \cup \{v\})$.
So $N(u) \cup N(v) = \{u,v,v_1,\ldots,v_n\}$ with $\{v_1,\ldots,v_n\} = N(v)\backslash\{u\}$.
The desired formula is a result of combining Theorem \ref{prop: ekf} with Lemma
\ref{betti J cap K} and using the fact that $\beta_{i,j}((uv)) = 0$ if $i \geq 1$.
\end{proof}

\begin{corollary} \label{reg pdim 2}
With the hypotheses and notation as in Theorem \ref{recursive formula 3}, we have
\begin{enumerate}
\item[$(i)$] $\reg(\I(G)) = \max\{2,\reg(\I(G \backslash e)), \reg(\I(H))+1\}.$
\item[$(ii)$] $\pdim(\I(G)) = \max\{\pdim(\I(G \backslash e)),\pdim(\I(H))+n+1\}$.
\end{enumerate}
\end{corollary}
\begin{proof}  Set $L = J \cap K$ and $G' = G \backslash e$.
By Corollary \ref{reg pdim} we have
\[\reg(\I(G)) =
\max\{\reg((uv)),\reg(\I(G')), \reg(L) -1 \}.\]
Since $\reg((uv)) = 2$,  we only need to
verify that $\reg(L) = \reg(\I(H)) + 2$. This is indeed true by 
Lemma \ref{betti J cap K}.  This proves $(i)$.
Similarly,  Corollary \ref{reg pdim} implies
\[\pdim(\I(G)) = \max\{\pdim((uv)), \pdim(\I(G')), \pdim(L) + 1 \}.\]
Now clearly $\pdim((uv)) = 0$.  By Lemma \ref{betti J cap K} we have
\[\pdim(L) = \pdim(R/((v_1,\ldots,v_n) + \I(H)))-1
= n + \pdim(R/\I(H)) - 1.\]
Since $\pdim(R/\I(H)) = \pdim(\I(H))+1$, the assertion $(ii)$ follows.
\end{proof}

\begin{example}
The above corollary implies 
that removing a splitting edge $e$ may decrease both the regularity and projective dimension,
that is, $\reg(\I(G)) \geq \reg(\I(G\backslash e))$ and $\pdim(\I(G)) \geq
\pdim(\I(G\backslash e)).$ However, if $e$ is not a splitting edge,
then it may happen  that $\reg(\I(G\backslash e))$, respectively
$\pdim(\I(G\backslash e))$, is larger than $\reg(\I(G))$, respectively $\pdim(\I(G))$.
For example, consider the graph $G$ below:
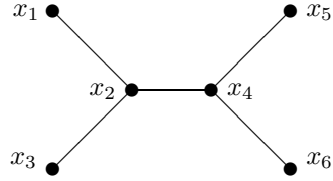
\[
\begin{picture}(100,80)(0,-35)
\put(0,0){\circle*{5}}
\put(-16,-2){$x_2$}
\put(30,0){\circle*{5}}
\put(36,-2){$x_4$}
\put(60,30){\circle*{5}}
\put(66,28){$x_5$}
\put(60,-30){\circle*{5}}
\put(66,-28){$x_6$}
\put(-30,-30){\circle*{5}}
\put(-46,-28){$x_3$}
\put(0,0){\line(1,0){30}}
\put(-30,30){\circle*{5}}
\put(-45,28){$x_1$}
\put(-30,30){\line(1,-1){30}}
\put(-30,-30){\line(1,1){30}}
\put(30,0){\line(1,1){30}}
\put(30,0){\line(1,-1){30}}
\end{picture}\]
The edge $x_2x_4$ is not a splitting edge. The resolution of $\I(G)$ is 
\[0 \rightarrow R^2(-4) \rightarrow R^6(-3) \rightarrow R^5(-2) \rightarrow \I(G)
\rightarrow 0\]
and the resolution of $\I(G\backslash e)$ is 
\[0 \rightarrow R(-6) \rightarrow R^4(-5) \rightarrow R^2(-3) \oplus R^4(-4)
\rightarrow R^4(-2) \rightarrow \I(G\backslash e) \rightarrow 0.\]
We have $\pdim(\I(G\backslash e)) = 3 > 2 = \pdim(\I(G))$
and $\reg(\I(G\backslash e)) = 3 > 2 = \reg(\I(G))$.
\end{example}

We end this section by using Theorem \ref{recursive formula 3}
to give new proofs for known results about the
the graded Betti numbers of forests.  We begin be recovering the
recursive formula of \cite{J,JK} found via a different means.
In fact, our result is more general since it applies to any leaf
of $G$, while  \cite{J,JK} required that  a special leaf  be removed.

\begin{corollary}\label{treeformula}
Let $e=uv$ be any leaf of a forest $G$.
If $\deg v =n$ and $N(v) = \{u,v_1,\ldots,v_{n-1}\}$, then for $i \geq 1$ and $j \geq 0$
\[\beta_{i,j}(\I(G)) = \beta_{i,j}(\I(T))
+ \sum_{l=0}^i \binom{n-1}{l} \beta_{i-l-1,j-2-l}(\I(H))\]
where $T = G \backslash e = G \backslash \{u\}$ and $H = G \backslash
\{u,v,v_1,\ldots,v_{n-1}\}$.
\end{corollary}

\begin{proof}
The hypotheses imply that $\deg u =1$.  Since $N(u) \subseteq (N(v) \cup \{v\})$,
$uv$ is a splitting edge.  Now
apply Theorem \ref{recursive formula 3}.
\end{proof}


Applying Corollary \ref{reg pdim 2} allows us to rediscover Theorem 4.8 of \cite{JK}.
\begin{corollary}\label{tree pdim}
With the notation as in the previous corollary,
\[\pdim(\I(G)) = \max\{\pdim(\I(T)),\pdim(\I(H))+n\}.\]
\end{corollary}

We say two edges $u_1v_1$ and $u_2v_2$ of a simple graph $G$ are {\it disconnected} if
(a) $\{u_1,v_1\} \cap \{u_2,v_2\} = \emptyset$, and
(b) $u_1u_2,u_1v_2,v_1u_2,v_1v_2$ are not
edges of $G$.  When $G$ is a forest, Theorem \ref{recursive formula 3} 
can be used to give a new proof of Zheng's result  (Theorem 2.18 of \cite{Z}) 
relating  $\reg(\I(G))$ to the number of disconnected edges.

\begin{corollary} \label{tree reg}
Let $G$ be a forest with edge ideal $\I(G)$.  Then
$\reg(\I(G)) = j + 1$  where $j$ is the maximal number of pairwise disconnected edges in $G$.
\end{corollary}

\begin{proof} We use induction on $|E_G|$.  The formula is clearly true for $|E_G| = 1$.

Suppose $|E_G| > 1$, and let $e = uv$ be any leaf of $G$ with $\deg u = 1$.  By
Corollary \ref{reg pdim 2} we have
\[\reg(\I(G)) = \max\{2,\reg(\I(T)),\reg(\I(H))+1\}\]
where $T = G \backslash e = G \backslash \{u\}$ and $H = G \backslash (\{v\} \cup N(v))$.
By induction $\reg(\I(T)) = j_1+1$ where $j_1$ is the maximal number
of pairwise disconnected edges of $T$, and $\reg(\I(H)) = j_2+1$
where $j_2$ is the maximal number of pairwise disconnected edges of $H.$
Since $\I(T)$ has at least one edge, $j_1 + 1 \geq 2$.  So
\[\reg(\I(G)) = \max\{j_1+1,j_2+2\}.\]
If we let $j$ denote the maximal number of pairwise disconnected edges of $G$,
then to complete the proof it suffices for us to show that $j = \max\{j_1,j_2+1\}$.

Let $\mathcal{E}_1$ be the set of the $j_1$ pairwise disconnected edges of $T$.  The
edges of $\mathcal{E}_1$ are also a set of pairwise disconnected edges of $G$.
Thus $|\mathcal{E}_1| = j_1 \leq j$.  If $\mathcal{E}_2$
is a set of $j_2$ pairwise disconnected edges of $H$, we claim that
$\mathcal{E}_2 \cup \{uv\}$ is a set of pairwise disconnected edges of $G$.
Indeed, $uv$ does not share a vertex with any edge in $H$. The only edges
that are adjacent to $uv$ are $vv_i$ with $v_i \in N(v)\backslash\{u\}$.
No edge of $\mathcal{E}_2$ can share a vertex with these edges since none
of the vertices of $N(v)$ belong to $H$.  Thus $|\mathcal{E}_2 \cup \{uv\}|
= j_2 + 1 \leq j$. Thus $j \geq \max\{j_1,j_2+1\}$.

Suppose that  $j > \max\{j_1,j_2+1\}$.  Let $\mathcal{E}$ be a set of $j$
pairwise disconnected edges of $G$.  If $uv \not\in\mathcal{E}$,
then $\mathcal{E}$ is also a set of pairwise disconnected edges of $T$,
and so $j = |\mathcal{E}| \leq j_1$, a contradiction.  If $uv \in \mathcal{E}$,
then $\mathcal{E} \backslash \{uv\}$ is a set of  pairwise disconnected
edges of $H$.  But this would imply that $j-1 \leq j_2$, again a
contradiction.  Hence
$j = \max\{j_1,j_2+1\}$.
\end{proof}


\section{Splitting Vertices} \label{sec-vertex}

Let $G$ be a simple graph, and let $v$ be a vertex of $G$
with $N(v) = \{v_1,\ldots,v_d\}$.  This
section complements the results of the previous section by determining
when $\I(G) = J + K$ with $J =  (vv_1,\ldots,vv_d)$ and $K = \I(G\backslash\{v\})$
is a splitting of $\I(G)$.

If $v$ is an isolated vertex of $G$, then  $\beta_{i,j}(\I(G))
= \beta_{i,j}(\I(G \backslash \{v\}))$ for all $i,j\geq 0$.  
If $\deg v = d > 0$
and if $G \backslash\{v\}$ consists of isolated vertices,
then $G = \K_{1,d}$, the complete bipartite graph of
size $1,d$; in this case the graded Betti numbers of $\I(G)$ follow
from Theorem \ref{resk1d}.  If $v \in V_G$ is in neither 
of these two cases,
we give it the following name.

\begin{definition} A vertex $v \in V_G$ is a {\it splitting vertex}
if $\deg v = d > 0$ and $G \backslash \{v\}$ is not the graph
of isolated vertices.
\end{definition}

This name makes sense in light of the following theorem.

\begin{theorem}\label{prop: edgeidealssplit}
Let $v$ be a splitting vertex of $G$ with $N(v) = \{v_1,\ldots,v_d\}$, 
and set $J = (vv_1,\ldots, vv_d)$ and $K = \I(G\backslash\{v\})$.
Then $\I(G) = J+K$ is a splitting of $\I(G)$.
\end{theorem}

\begin{proof}
It is clear that  $\I(G) = J + K$.  As well, 
$\G(\I(G)) = \G(J) \cup \G(K)$ is a disjoint union
because $v$ divides all elements of $\G(J)$ but divides no element of $\G(K)$.

Now consider the ideal $J \cap K = (vv_1,\ldots,vv_d) \cap \I(H)$
where $H = G \backslash \{v\}$.  Then
\[J \cap K = ( \{\lcm(m_1,m_2) ~|~
m_1 \in \{vv_1,\ldots,vv_d\}, ~m_2 \in \G(\I(H))\}).\]
Thus
\begin{eqnarray}
\G(J \cap K) & = & \{ vv_iv_j ~|~ v_iv_j \in E_G\} \cup
\{vv_iy_j ~|~ v_iy_j \in E_G\} \cup \nonumber \\
&& \{vv_iy_jy_k ~|~ y_jy_k \in E_G ~\mbox{but}~ v_iy_j,
 v_iy_k \not\in E_G\} \label{JKdescription}
\end{eqnarray}
where $y_i$ denotes a vertex in $V_G \backslash \{v,v_1,\ldots,v_d\}$.
Note that the three sets are disjoint.

We define a splitting function $\G(J\cap K) \rightarrow \G(J) \times \G(K)$
as follows.  If $w \in \G(J\cap K)$, then define $\phi: \G(J\cap K)
\rightarrow \G(J)$ and $\psi:\G(J\cap K) \rightarrow \G(K)$ by
\[
\phi(w) =
\left\{
\begin{array}{ll}
vv_i & \mbox{if $w = vv_iv_j$ and $i < j$} \\
vv_i & \mbox{if $w = vv_iy_j$}\\
vv_i & \mbox{if $w = vv_iy_jy_k$}
\end{array}
\right.
\mbox{ and }
\psi(w) =
\left\{
\begin{array}{ll}
v_iv_j & \mbox{if $w = vv_iv_j$}\\
v_iy_j & \mbox{if $w = vv_iy_j$}\\
y_jy_k & \mbox{if $w = vv_iy_jy_k$.}
\end{array}
\right.
\]
By construction, the map given by $w \mapsto (\phi(w),\psi(w))$
has the property that $w = \lcm(\phi(w),\psi(w))$.
It suffices to verify condition $(b)$ of $(2)$ in Definition \ref{defn: split}.

So, suppose $S \subseteq \G(J \cap K)$. If $S$ contains a monomial divisible by some 
variable $y \not\in \{v, v_1, \dots, v_d\}$, then $\lcm(\phi(S))$ strictly divides 
$\lcm(S)$ since $y$ does not divide $\lcm(\phi(S))$. Otherwise, we must have 
$S \subseteq \{ vv_iv_j ~|~ v_iv_j \in E_G\}$. In this case, let $f$ be the maximal 
index such that $v_f$ appears in a monomial of $S$. Then, by the definition of $\phi$, 
$v_f$ does not divide $\phi(w)$ for any $w \in S$. Thus, $v_f$ does not divide 
$\lcm(\phi(S))$. Therefore, $\lcm(\phi(S))$ strictly divides $\lcm(S)$.
It is clear that $\lcm(\psi(S))$ strictly divides $\lcm(S)$
because $v$ does not divide $\lcm(\psi(S))$.  The theorem is proved.
\end{proof}

The following result is an immediate consequence of our description in (\ref{JKdescription}).

\begin{corollary} \label{betti intersection 2}
With the notation as in the previous theorem, set
\begin{eqnarray*}
G_i & := & G \backslash (N(v) \cup N(v_i)) ~\mbox{for $i = 1,\ldots,d,$ and } \\
G_{(v)} & := & G_{\{v_1,\ldots,v_d\}} \cup 
\{e\in E_G ~|~ \mbox{$e$ incident to one of $v_1,\ldots,v_d$,
but not $v$}\}.
\end{eqnarray*}
Then
\[J \cap K = v\I(G_{(v)}) + vv_1\I(G_1) + vv_2\I(G_2)+\cdots + vv_d\I(G_d).\]
\end{corollary}

Theorem \ref{prop: edgeidealssplit} gives us
some partial results on how the projective dimension and regularity
behave under removing any (splitting or non-splitting) vertex.

\begin{corollary} \label{reg pdim 3}
Let $G$ be a simple graph, and let $v \in V_G$ be any vertex.  Then
\begin{enumerate}
\item[$(i)$] $\reg(\I(G)) \geq \max\{2,\reg(\I(G\backslash\{v\}))\}.$
\item[$(ii)$]$\pdim(\I(G)) \geq \max\{d-1,\pdim(\I(G\backslash\{v\}))\}$ where
$d = \deg v$.
\end{enumerate}
\end{corollary}

\begin{proof}
If $v$ is not a splitting vertex, then $(i)$ and $(ii)$ are immediate from
the fact that either $\I(G) = \I(G\backslash \{v\})$, or $\I(G) = \I(\K_{1,d})$
and $\I(G\backslash\{v\}) = (0)$.
If $v$ is a splitting vertex, then
$\I(G) = (vv_1,\ldots,vv_d) + \I(G\backslash \{v\})$ is a splitting.
Now use Corollary \ref{reg pdim} and the fact that
$\reg(\I(\K_{1,d})) =2$ and $\pdim(\I(\K_{1,d})) = d-1$.
\end{proof}

\begin{remark}
Jacques proved (Proposition 2.1.4 of \cite{J}) statement $(ii)$
when the vertex $v$ is a terminal vertex, i.e.,
adjacent to at most one other vertex of $G$.
\end{remark}

Applying Theorems \ref{prop: ekf} and \ref{prop: edgeidealssplit}
and Corollary \ref{betti intersection 2} we obtain our next main result.

\begin{theorem}\label{recursive formula 2}
Let $v$ be a splitting vertex of $G$ with $N(v) = \{v_1,\ldots,v_d\}$. Let $G_{(v)}$ and $G_i$ ($i=1, \dots, d$) be defined as in Corollary \ref{betti intersection 2}. Then
\[\beta_{i,j}(\I(G)) = \beta_{i,j}(\I(\K_{1,d})) + \beta_{i,j}(\I(G \backslash \{v\})) +
\beta_{i-1,j}(L)\]
where  $L = v\I(G_{(v)}) + vv_1\I(G_1) + \cdots +vv_d\I(G_d)$ and $\K_{1,d}$ 
is the complete bipartite graph of size $1,d$.
\end{theorem}

Theorem \ref{recursive formula 2} allows us to give a new combinatorial
proof for an interesting result due to Eisenbud, et al. \cite{eghp}. We say that a
cycle $C = (x_1x_2 \dots x_qx_1)$ of $G$ has a {\it chord} if there exists some $
j \not\equiv i+1 (\text{mod} \ q)$ such that $x_ix_j$ is an edge of $C$. We call a
cycle $C$ a {\it minimal cycle} if $C$ has length at least 4 and contains no chord. An
ideal $I$ is said to satisfy {\it property $N_{2,p}$} for some $p \ge 1$ if $I$ is
generated by quadratics and its minimal free resolution is linear up to the $p$th step,
i.e., $\beta_{i,j}(I) = 0$ for all $0 \le i < p$ and $j > i+2$.

\begin{corollary}[see Theorem 2.1 of \cite{eghp}] \label{np property}
Let $G$ be a simple graph with edge ideal $\I(G)$. Then $\I(G)$ satisfies property $N_{2,p}$ with $p > 1$
if and only if every minimal cycle in $G^c$ has length $\ge p+3$.
\end{corollary}

\begin{proof} We use induction on $n = |V_G|$. Our assertion is vacuously true for $n \le 3$.
Suppose $n \ge 4$. We may assume that $G$ has no
isolated vertices. Since the edge ideal of the complete bipartite graph $\K_{1,n-1}$ has a
linear resolution by Theorem \ref{resk1d}, our statement is also vacuously true
in this case.

Suppose $G$ is
not the complete bipartite graph $\K_{1,n-1}$. Clearly $G$ now has a splitting vertex, say
$v$. Set $N(v) = \{v_1, \dots, v_d\}$, and let $G_i = G \backslash (N(v) \cup N(v_i))$ for $i=1,
\dots, d$, and $G_{(v)} = G_{\{v_1, \dots, v_d\}} \cup 
\{ e \in E_G ~|~ \mbox{$e$ is incident to one of $v_1,\ldots,v_d$ but not $v$}\}$. By Theorem \ref{recursive formula 2}
(and Corollary \ref{betti intersection 2}) we have that
\begin{align}
\beta_{i,j}(\I(G)) = \beta_{i,j}(J)+\beta_{i,j}(K)+\beta_{i-1,j}(L) \label{np equation}
\end{align}
where $J = (vv_1, \dots, vv_d)$, $K = \I(G \backslash \{v\})$ and $L = J \cap K = v\I(G_{(v)}) + \sum_{i=1}^d vv_i\I(G_i)$.

It follows from (\ref{np equation}) that $\I(G)$ satisfies property $N_{2,p}$ if and only if $J$
and $K$ satisfy property $N_{2,p}$, and $L$ satisfies property $N_{3,p-1}$. Observe further that
$L$ satisfies property $N_{3,p-1}$ if and only if $L = v\I(G_{(v)})$ and $\I(G_{(v)})$ satisfies
 property $N_{2,p-1}$. Since $J$ has a linear minimal free resolution, $J$ always satisfies
property $N_{2,p}$. By the induction hypothesis, $K$ satisfies property $N_{2,p}$ if and only if
 every minimal cycle of $(G \backslash \{v\})^c$ has length $\ge p+3$. It can be seen that
$(G \backslash \{v\})^c = G^c \backslash \{v\}$. Thus, it remains to prove that $L = v\I(G_{(v)})$
and $\I(G_{(v)})$ satisfies property $N_{2,p-1}$ if and only if every minimal cycle of $G^c$
containing $v$ has length $\ge p+3$.

Suppose first that $L = v\I(G_{(v)})$ and $\I(G_{(v)})$ satisfies property $N_{2,p-1}$.
Consider $C = (vx_1 \dots x_lv)$ an arbitrary minimal cycle in $G^c$ containing $v$
(and thus, $l \geq 3$). We shall
show that $C$ has length $\ge p+3$. Since $C$ is a minimal cycle, we have $vx_2, vx_3, \dots,
 vx_{l-1} \not\in G^c$. This implies that $vx_2, \dots, vx_{l-1} \in G$. Thus, $\{x_2, \dots,
x_{l-1}\} \subseteq \{v_1, \dots, v_d\}$. Also, since $vx_1, vx_l \in G^c$, we have $x_1, x_l
\not\in \{v_1, \dots, v_d\}$. This implies that $x_1x_l \not\in G_{(v)}$. Therefore, $x_1, \dots,
 x_l$ form either a minimal cycle or a triangle in $G_{(v)}^c$. Consider the case when $l \ge 4$.
 If $p=2$, then clearly $C$ has length $\ge p+3$. If $p > 2$ then by the induction hypothesis, since
$G_{(v)}$ does not contain $v$ and $\I(G_{(v)}^c)$ satisfies property $N_{2,p-1}$, every minimal
 cycle in $G_{(v)}^c$ must have length $\ge p+2$. Hence, $l \ge p+2$, whence $C$ has length
$\ge p+3$. It remains to consider the case when $l=3$. Since $C$ is a minimal cycle, $x_1x_3$
is not a chord of $C$. This means that $x_1x_3 \in G$. Furthermore, as shown, $x_1, x_3 \not\in
 \{v_1, \dots, v_d\}$ and $x_2x_1, x_2x_3 \in G^c$. This implies that $vx_2x_1x_3 \in J \cap K = L$
and $vx_2x_1x_3 \not\in v\I(G_{(v)})$, a contradiction to the fact that $L = v\I(G_{(v)})$.

Conversely, suppose that every minimal cycle of $G^c$ containing $v$ has length $\ge p+3$.
We need to prove: (a) $L = v\I(G_{(v)})$, and (b) $\I(G_{(v)})$ has property $N_{2,p-1}$.

To prove (a) we observe that if $v\I(G_{(v)}) \subsetneq L$ then there exists an edge
$e = uw \in G$ such that $u,w \not\in \{v, v_1, \dots, v_d\}$.  But then $(vuv_iw)$
forms a minimal cycle of length 4 in $G^c$, contradicting the assumption that every
 minimal cycle of $G^c$ has length $\ge p+3$ for $p > 1$.

To prove (b) we observe that $\I(G_{(v)})$ is generated by quadratics, so (b) is true for $p=2$.
Assume that $p > 2$. By induction, we only need to show that every minimal cycle of $G_{(v)}^c$ has
 length $\ge p+2$. Consider an arbitrary minimal cycle $D = (x_1x_2 \dots x_lx_1)$ in $G_{(v)}^c$.
 Let $\{w_1, \dots, w_s\}$ be the set of vertices of $G \backslash \{v, v_1, \dots, v_d\}$ which
are adjacent to at least one of the vertices $\{v_1, \dots, v_d\}$. If there exist $1 \le i \not= j
\le s$ such that $w_i, w_j \in \{x_1, \dots, x_l\}$, then since $w_iw_j \not\in G_{(v)}$ (by
definition of $G_{(v)}$), $w_iw_j$ must be an edge of $D$ (otherwise $D$ would have a chord).
Without loss of generality, suppose $w_i = x_1$ and $w_j = x_l$. In this case, $(x_1x_2 \dots x_l
 v x_1)$ is a minimal cycle of $G^c$, which implies that $l+1 \ge p+3$, i.e., $l \ge p+2$. If
 there is at most one of $\{w_1, \dots, w_s\}$ belonging to $\{v_1, \dots, v_d\}$, then it is easy
to see that $D$ is a minimal cycle in $G^c$. This implies that $l \ge p+3 > p+2$. The result is
proved.
\end{proof}

Fr\"oberg's \cite{Fr}
main theorem is a special case of Corollary \ref{np property}.
 Recall that $\I(G)$ has a {\it linear resolution} if $\beta_{i,j}(\I(G)) = 0$
for all $j \neq i+2$. We say that $G$ is {\it chordal} if every cycle of length $> 3$ has a chord;
in other words, $G$ has no minimal cycles.

\begin{corollary}[see \cite{Fr}] \label{linear resolution} Let $G$ be a graph with edge ideal $\I(G)$.
Then $\I(G)$ has a linear resolution if and only if $G^c$ is a
chordal graph.
\end{corollary}

As a consequence of Theorem \ref{recursive formula 2}, we can give a new proof for
the formula of \cite {RVT} 
for the linear strand of $\I(G)$ when $G$ contains no minimal cycle of length 4 .

\begin{corollary}[see \cite{RVT}] \label{linear strand}
Let $G$ be a graph with no minimal cycle of length 4. Let $k_{i+2}(G)$ denote the number of $(i+2)$-cliques in $G$. Then, for any $i \geq 0$,
$$\beta_{i,i+2}(\I(G)) = \sum_{u \in V_G} {\deg u \choose i+1} - k_{i+2}(G).$$
\end{corollary}

\begin{proof} We have
$\sum_{u \in V_G} \deg u = 2 |E_G| = 2k_2(G).$
Thus, our statement is true for $i = 0$. Assume that $i \ge 1$. We shall use induction on $n = |V_G|$.

It is easy to verify the statement for $n \le 3$. Assume that $n \ge 4$. The
statement for the complete bipartite graph $\K_{1,n-1}$ follows
by Theorem \ref{resk1d}, so we may assume that $G$ is not the complete
 bipartite graph $\K_{1,n-1}$. This guarantees that $G$ contains a splitting vertex, say $v$. As
before, let $N(v) = \{v_1, \dots, v_d\}$, let $G_i = G \backslash (N(v) \cup N(v_i))$ for $i=1,
\dots, d$, and let $G_{(v)} = G_{\{v_1, \dots, v_d\}} \cup \{ e \in E_G \ | \ e \ \text{incident to}
 \ v_i \ \text{for some} \ i=1, \dots, d \ \text{but not}\ v\}$. By Theorem \ref{recursive formula 2} (and Corollary
 \ref{betti intersection 2}), we have
$\beta_{i,i+2}(\I(G)) = \beta_{i,i+2}(J) + \beta_{i,i+2}(K) + \beta_{i-1,i+2}(L)$
where $J = (vv_1, \dots, vv_d)$, $K = \I(G \backslash \{v\})$, and $L = v\I(G_{(v)}) + \sum_{i=1}^d
 vv_i\I(G_i)$.  For each $i=1, \dots, d$, the ideal $vv_i\I(G_i)$ is generated by monomials of
degree 4.  Thus
 the linear strand of $L$ is the same as that of $v\I(G_{(v)})$ (or equivalently, that of
$\I(G_{(v)})$). Thus, we have
\begin{align}
\beta_{i,i+2}(\I(G)) = \beta_{i,i+2}(J)+\beta_{i,i+2}(K)+\beta_{i-1,i+1}(\I(G_{(v)})).
\label{linear equation}
\end{align}

For simplicity, let $G' = G \backslash \{v\}$ and $G'' = G_{(v)}$. Let $W = \{w_1, \dots, w_s\}$ be
the set of vertices of $G \backslash \{v, v_1, \dots, v_d\}$ which are adjacent to at least one
of the vertices of $N(v) = \{v_1, \dots, v_d\}$. Let $\{v_{j1}, \dots, v_{jl_j}\}$ be the set of
 vertices of $\{v_1, \dots, v_d\}$ that are adjacent to $w_j$, for $j=1, \dots, s$. By the induction
hypothesis we have
\begin{align}
\footnotesize
\beta_{i-1,i+1}(\I(G'')) & = \sum_{u \in N(v)} {\deg_{G''} u \choose i} + \sum_{u \in W} {\deg_{G''} u \choose i} - k_{i+1}(G'') \nonumber \\
& = \sum_{u \in N(v)} {\deg_G u - 1 \choose i} + \sum_{j=1}^s {l_j \choose i} - k_{i+1}'(G'') -
k_{i+1}''(G'') \label{betti:G''}
\normalsize
\end{align}
where $k_{i+1}'(G'')$ denotes the number of $(i+1)$-cliques $G''$ not containing any of the
vertices in $W$ and $k_{i+1}''(G'')$ denotes the number of $(i+1)$-cliques of $G''$ containing at
least one vertex in $W$. Observe that for each $j=1, \dots, s$, since $G$ contains no minimal cycle of length 4
and $vw_j \not\in E_G$, we must have $v_{jt_1}v_{jt_2} \in E_G$ for any $1 \le t_1 \not= t_2 \le l_j$.
 This implies that $G_{\{v_{j1}, \dots, v_{jl_j}\}}$ is the complete graph on $l_j$ vertices for
 any $j=1, \dots, s$. Moreover, $w_lw_t \not\in E_{G''}$ for any $1 \le l \not= t \le s$.
Therefore, each $(i+1)$-clique of $G''$ containing some vertices in $W$ contains exactly one.
We must have
$k_{i+1}''(G'') = \sum_{j=1}^s {l_j \choose i}.$
This, together with (\ref{betti:G''}), gives
\begin{align}
\beta_{i-1,i+1}(\I(G'')) 
= \sum_{u \in N(v)} {\deg_G u - 1 \choose i} - k_{i+2}(G_{\{v, v_1, \dots, v_d\}}). \label{k'}
\end{align}
By induction we also have
\begin{align}
\beta_{i,i+2}(K) & = \sum_{u \in V_G \backslash \{v, v_1, \dots, v_d\}} {\deg_G u \choose i+1} + \sum_{u \in N(v)} {\deg_G u - 1 \choose i+1} - k_{i+2}(G'). \label{betti:K}
\end{align}
It can further be seen that
\begin{align}
\beta_{i,i+2}(J) = {\deg_G v \choose i+1}. \label{betti:J}
\end{align}

Now (\ref{linear equation}), (\ref{k'}), (\ref{betti:K}), and (\ref{betti:J}) combine to give us
\begin{align*}
\beta_{i,i+2}(\I(G)) = \sum_{u \in V_G} {\deg_G u \choose i+1} - k_{i+2}(G') - k_{i+2}(G_{\{v,v_1, \dots, v_d\}}).
\end{align*}
Observe that an $(i+2)$-clique in $G$ either contains $v$ (so it is an $(i+2)$-clique of
$G_{\{v,v_1,\dots,v_d\}}$) or is a $(i+2)$-clique of $G' = G \backslash \{v\}$. Hence,
$$k_{i+2}(G) = k_{i+2}(G_{\{v,v_1,\dots,v_d\}}) + k_{i+2}(G')$$
and thus the result is proved.
\end{proof}


\section{Facet Ideals and Splitting Facets}

In this section we extend our method to the study of arbitrary square-free monomial ideals.
Let $\Delta$ be a simplicial complex on the vertex set $V_\Delta = \{x_1,\ldots,x_n\}$.
Let $\I(\Delta)$ be the
facet ideal of $\Delta$ in  $R = k[x_1,\ldots,x_n]$. Recall that, by
abuse of notation, we will use $F$ to denote a facet of $\Delta$ 
and the monomial $\prod_{x \in F} x$ in $\I(\Delta)$.

\begin{definition} \label{3conn}
Let $F$ be a facet of $\Delta$.
The {\it connected component} of $F$ in $\Delta$, denoted  $\lk_\Delta(F)$,
is defined to be the connected component of $\Delta$ containing $F$.
If $\lk_\Delta(F) \backslash F = \langle G_1, \dots, G_p \rangle$, then we define
the {\it reduced connected component} of $F$ in $\Delta$, denoted by $\rlk_{\Delta}(F)$, to
be the simplicial complex whose facets are given by
$G_1 \backslash F, \dots, G_p \backslash F$, where if there exist $G_i$ and $G_j$ such
that $\emptyset \not= G_i \backslash F \subseteq G_j \backslash F$, then we shall disregard
the bigger facet $G_j \backslash F$ in $\rlk_\Delta(F)$.
\end{definition}

Let $F$ be a facet of $\Delta$. Let $\Delta' = \Delta \backslash F$ be the simplicial
complex whose facet set is $\F(\Delta) \backslash F$. Let
$J = (F)$ and $K = \I(\Delta').$
Note that $\G(\I(\Delta))$ is the disjoint union of $\G(J)$ and $\G(K)$. We are interested in
finding $F$ such that $\I(\Delta) = J + K$ gives a splitting for $\I(\Delta)$.

\begin{definition} \label{defn: splitting facet}
With the above notation, we shall call $F$ a {\it splitting facet} of $\Delta$ if $\I(\Delta) = J + K$
is a splitting of $\I(\Delta)$.
\end{definition}

\begin{lemma} \label{3generators}
With the above notation, we have
$$J \cap K = (F)(\I(\rlk_\Delta(F)) + \I(\Omega))$$
where $\Omega$ denotes the simplicial complex $\Delta \backslash \lk_\Delta(F).$
\end{lemma}

\begin{proof} Since both $J$ and $K$ are monomial ideals, we have
$$J \cap K = (\{ \lcm(F,G) \ | \ G \in \G(K)\}).$$
It is easy to see that if $H$ is a facet of $\rlk_\Delta(F)$ and if $G$ is a facet of
$\lk_\Delta(F)$ such that $G \backslash F = H$, then $G \not= F$ and
$FH = \lcm(F,G) \in J \cap K$. Thus, $(F)\I(\rlk_\Delta(F)) \subseteq J \cap K$.
Also, for any facet $H$ of $\Omega$, $FH = \lcm(F,H) \in J \cap K$. Thus,
$(F)\I(\Omega) \subseteq J \cap K$, and hence
$$ (F)(\I(\rlk_\Delta(F)) + \I(\Omega)) \subseteq J \cap K.$$

For the other inclusion, note that each $G \in \G(K)$ corresponds to
a facet $G$ of $\Delta'$. There are two possibilities for this facet: (1) $G \in \lk_\Delta(F)$,
or (2) $G \not\in \lk_\Delta(F)$. It now follows from the construction of $\rlk_\Delta(F)$ and 
$\Omega$ that case (1) leads to 
$\lcm(F,G) \in (F)\I(\rlk_\Delta(F))$ and case (2) results in $\lcm(F,G) \in (F)\I(\Omega)$.
\end{proof}

\begin{lemma} \label{3intBetti}
With the same notation as in Lemma \ref{3generators}, for all $i \geq 1$ and all $j \geq 0$
$$\beta_{i-1,j}(J\cap K) = \sum_{l_1=0}^i \sum_{l_2=0}^{j-|F|} \beta_{l_1-1,l_2}(\I(\rlk_\Delta(F)))
\beta_{i-l_1-1,j-|F|-l_2}(\I(\Omega)).$$
\end{lemma}

\begin{proof} Let $L = \I(\rlk_\Delta(F)) + \I(\Omega).$ It follows from Lemma
\ref{3generators} that $J \cap K = FL$. Since none of the variables in $F$ are present in
$L$, $F$ is not a zero-divisor of $R/L$. As a consequence, we have for all $i \geq 1$
\begin{align*}
\beta_{i-1,j}(FL) = \beta_{i-1,j-|F|}(L) = \beta_{i,j-|F|}(R/L).
\end{align*}

Now we notice that $\rlk_\Delta(F)$ and $\Omega$ do not share any common vertices. 
The statement, therefore, follows by applying Remark \ref{rmk-tensor}.
\end{proof}

The following result gives a recursive like formula for the graded Betti numbers
of the facet ideal of a simplicial complex in terms of the Betti numbers
of facet ideals of subcomplexes.  This
result is a higher dimension analogue of Theorem \ref{recursive formula 3}.

\begin{theorem} \label{simplicial:betti}
Let $F$ be a splitting facet of a simplicial complex $\Delta$. Let $\Delta' = \Delta \backslash F$ 
and $\Omega = \Delta \backslash \lk_\Delta(F).$ Then, for all $i \geq 1$ and $j \geq 0$,
$$\beta_{i,j}(\I(\Delta)) = \beta_{i,j}(\I(\Delta')) + \sum_{l_1=0}^i
\sum_{l_2=0}^{j-|F|} \beta_{l_1-1,l_2}(\I(\rlk_\Delta(F))) \beta_{i-l_1-1,j-|F|-l_2}(\I(\Omega)).$$
\end{theorem}

\begin{proof} By definition, $\I(\Delta) = J + K$ is a splitting of $\I(\Delta)$. 
The conclusion now follows
from Theorem \ref{prop: ekf} and Lemma \ref{3intBetti} and the fact that $\beta_{i,j}(J) = 0$
if $i \geq 1$.
\end{proof}

We will now show that our formula in Theorem \ref{simplicial:betti} is recursive when $G$ is a forest.
To do so, we first show that a leaf of $\Delta$ is a splitting facet.
Recall that if $F$ is a leaf of $\Delta$, then $F$ must have a vertex that does not belong to
any other facet of the simplicial complex (see Remark 2.3 of \cite{faridi:2002}).

\begin{theorem} \label{3split}
If $F$ is a leaf of $\Delta$, then $F$ is a splitting facet of $\Delta$.
\end{theorem}

\begin{proof} We need to show that if $F$ is a leaf of $\Delta$,
 then $\I(\Delta) = J+K$ with $J = (F)$ and $K = \I(\Delta\backslash F)$
is a splitting of $\I(\Delta)$. Without loss of generality, we may assume
 that $F = \{ x_1, \dots, x_l \}$. We
shall construct a splitting function
$s: \G(J \cap K) \rightarrow \G(J) \times \G(K)$
for $\I(\Delta)$.

Suppose $L \in \G(J \cap K)$. Let $\M_L = \{ G \in \G(K) \ | \ \lcm (F,G) = L \}$. For each
 $G \in \M_L$, we order the elements of $G \cap F$ by the increasing order of their indexes and
 view $G \cap F$ as an ordered word of the alphabet $\{ x_1, \dots, x_l \}$. Let $G_L \in \M_L$
 be such that $G_L \cap F$ is minimal with respect to the lexicographic word ordering.
 Clearly, $G_L$ is uniquely determined by $L$. Our splitting function $s$ is defined as
follows. For each $L \in \G(J \cap K)$,
$$s(L) = (\phi(L), \psi(L)) = (F, G_L).$$
We need to verify that $s$ satisfies conditions (a) and (b) of Definition \ref{defn: split}.
Indeed, condition (a) follows obviously from the definition of the function $s$. Suppose
$S \subseteq \G(J \cap K)$. Since $F$ is a leaf of $\Delta$, there exists a vertex $u \in F$
such that $u$ is not in any other facet of $\Delta$. This implies that $u$ does not divide
 $\lcm(\psi(S))$. Yet, since $u$ is in $F$, $u$ divides $\lcm(S)$. Thus, $\lcm(\psi(S))$
strictly divides $\lcm(S)$. On the other hand, it is also clear that for any $G \in \G(K)$,
$F$ strictly divides $\lcm(F,G)$, so $\lcm(\phi(S)) = F$ strictly divides $\lcm(S)$. The
result is proved.
\end{proof}

Because $\lk_\Delta(F)$, $\Delta \backslash F$ and $\Delta \backslash \lk_\Delta(F)$ are
subcomplexes of $\Delta$, it follows directly from the definition that if $\Delta$ is a forest,
then so are $\lk_\Delta(F)$, $\Delta \backslash F$ and $\Delta \backslash \lk_\Delta(F)$.
Thus, to show that  our formula in Theorem \ref{simplicial:betti} is recursive when $G$ is a forest,
it suffices to show that $\rlk_\Delta(F)$ is also a forest.

\begin{lemma} \label{3forest}
Let $F$ be a facet of a forest $\Delta$. Then $\rlk_\Delta(F)$ is a forest.
\end{lemma}

\begin{proof} Suppose $\Xi = \langle G_1, \dots, G_l \rangle$ is a connected component of
$\rlk_\Delta(F)$, where $G_i = F_i \backslash F$ and $F_i$ is a facet of $\lk_\Delta(F)$ for
 all $i = 1, \dots, l$. We shall show that $\Xi$ has a leaf. Indeed, it is easy to see that
$\Theta = \langle F_1, \dots, F_l \rangle$ is a connected subcomplex of $\lk_\Delta(F)$. As
observed, since $\Delta$ is a forest, so is $\lk_\Delta(F)$. Thus, $\Theta$ has a leaf.
Without loss of generality, assume that $F_1$ is a leaf of $\Theta$. That is, either $l = 1$ or
there exists another facet of $\Theta$, say $F_2$, such that $F_1 \cap H \subseteq F_1 \cap F_2$
for any facet $H \not= F_1$ of $\Theta$. If $l = 1$, then clearly $G_1$ is a leaf of $\Xi$.
Suppose $l > 1$. It is easy to see that $G_1 \cap (H \backslash F) = (F_1 \backslash F) \cap
(H \backslash F) = (F_1 \cap H) \backslash F \subseteq (F_1 \cap F_2) \backslash F =
(F_1 \backslash F) \cap (F_2 \backslash F) = G_1 \cap G_2$. Thus, $G_1$ is also a leaf of
 $\Xi$. We have just shown that $\Xi$ has a leaf in any case. The lemma is proved.
\end{proof}

We can generalize Corollary \ref{treeformula} by giving a recursive formula for
simplicial trees.

\begin{theorem} \label{simplicial:recursive}
Let $\Delta$ be a simplicial forest.  For any leaf $F$ of $\Delta$,
$\Delta' = \Delta \backslash F$, $\Omega = \Delta \backslash \lk_\Delta(F)$,
and $\rlk_\Delta(F)$ are also simplicial forests.  Furthermore, the
numbers $\beta_{i,j}(\I(\Delta))$ for all $i \geq 1$ and $j \geq 0$ can be computed
recursively using the formula
$$\beta_{i,j}(\I(\Delta)) = \beta_{i,j}(\I(\Delta')) + \sum_{l_1=0}^i
\sum_{l_2=0}^{j-|F|} \beta_{l_1-1,l_2}(\I(\rlk_\Delta(F))) \beta_{i-l_1-1,j-|F|-l_2}(\I(\Omega)).$$
\end{theorem}

\begin{proof} 
Lemma \ref{3forest} and the discussion before this lemma imply the 
first statement.  The second statement follows from Theorems \ref{simplicial:betti}
and \ref{3split} because $\Delta'$, $\Omega$, and $\rlk_\Delta(F)$
all have must have a leaf, which implies their facet ideals can also be
split using Theorem \ref{simplicial:betti}.
\end{proof}

Theorem \ref{simplicial:recursive} can be used to find a nice formula for the linear strand of
facet ideals of pure forests, generalizing a result of Zheng \cite[Proposition 3.3]{Z} and Corollary \ref{linear strand}.
Recall that a simplicial complex $\Delta$ is said to be pure $(d-1)$-dimensional if
$\dim F = d-1$, i.e., $|F| = d$, for any facet $F$ of $\Delta$. For a face $G$ of
dimension $d-2$ of a pure $(d-1)$-dimensional simplicial complex $\Delta$ we define the
{\it degree} of $G$ in $\Delta$, written $\deg_\Delta(G)$, to be the cardinality of the set
 $\{ F \in \F(\Delta) \ | \ G \subseteq F \}$. Let $\A(\Delta)$ denote the set of
 $(d-2)$-dimensional faces of $\Delta$.

\begin{theorem} \label{3linear-strand}
Let $\Delta$ be a pure $(d-1)$-dimensional forest (for some $d \ge 2$). Then
$$\beta_{i,i+d}(\I(\Delta)) = \left\{ \begin{array}{lcl} |\F(\Delta)| & \text{if} & i = 0 \\
{\displaystyle \sum_{G \in \A(\Delta)} {\deg_\Delta(G) \choose i+1}} & \text{if} & i \ge 1.
\end{array} \right.$$
\end{theorem}

\begin{proof} The assertion is clear for $i = 0$. Suppose $i \ge 1$. Let $m = |\F(\Delta)|$
 be the number of facets of $\Delta$. We shall use induction on $m$. For $m = 1$, the assertion
is obviously true. Suppose that $m > 1$. Let $F$ be a leaf of $\Delta$, and
let $\Delta' = \Delta \backslash F$ and $\Omega = \Delta \backslash \lk_\Delta(F)$. By Theorem \ref{simplicial:recursive}, for
$i \ge 1$ we have
$$\beta_{i,i+d}(\I(\Delta)) = \beta_{i,i+d}(\I(\Delta'))+\sum_{l_1=0}^i \sum_{l_2=0}^i
\beta_{l_1-1, l_2}(\I(\rlk_\Delta(F))) \beta_{i-l_1-1,i-l_2}(\I(\Omega)).$$

Observe that since $d \ge 2$, $l_1-1 \ge l_1+1-d$. Thus, for any $l_2 = 0, \dots, i$, we
have either $l_2 \le l_1-1$ or $l_2 > l_1+1-d$. If $l_2 \le l_1-1$, then clearly
$\beta_{l_1-1, l_2}(\I(\rlk_\Delta(F))) = 0$. If $l_2 > l_1+1-d$, then $i-l_2 < (i-l_1-1)+d$. This and
the fact that $\Omega = \Delta \backslash \lk_\Delta(F)$ is also a pure $(d-1)$-dimensional
forest imply that $\beta_{i-l_1-1, i-l_2}(\I(\Omega)) = 0$ unless $l_1 = l_2 = i$ (in which case $\beta_{i-l_1-1,i-l_2}(\I(\Omega)) = \beta_{-1,0}(\I(\Omega)) = 1$).
Hence, we have
\begin{align}
\beta_{i,i+d}(\I(\Delta)) = \beta_{i,i+d}(\I(\Delta')) +
\beta_{i-1, i}(\I(\rlk_\Delta(F))). \label{3eq2}
\end{align}

Clearly, $\beta_{i-1,i}(\I(\rlk_\Delta(F)))$ forms the linear strand of $\I(\rlk_\Delta(F))$
and is given by ${\displaystyle s \choose i}$ for any $i \ge 1$, where $s$ is the number of
isolated vertices of $\rlk_\Delta(F)$. Since $F$ is a leaf of $\Delta$, there must exist a
 vertex $u \in F$ such that $u$ is not in any other facet of $\Delta$. Let
$H = F \backslash \{u\}$. Observe that $\{x\}$ (for some $x \not= u$) is an isolated vertex of
$\rlk_\Delta(F)$ if and only if $H \cup \{x\}$ is a facet of $\Delta$. This implies that
$s = \deg_\Delta(H) - 1$ (since $H = F \backslash \{u\}$ is not in $\rlk_\Delta(F)$). This,
together with the induction hypothesis, now gives
\allowdisplaybreaks

\begin{align*}
\beta_{i,i+d}(\I(\Delta)) & = \sum_{G \in \A(\Delta')} {\deg_{\Delta'}(G) \choose i+1} + {\deg_\Delta(H)-1 \choose i} \\
& = \sum_{G \in \A(\Delta) \backslash \{H\}} {\deg_\Delta(G) \choose i+1} + {\deg_\Delta(H)-1 \choose i+1} + {\deg_\Delta(H)-1 \choose i} \\
& = \sum_{G \in \A(\Delta)} {\deg_\Delta(G) \choose i+1}.
\end{align*}
The theorem is proved.
\end{proof}




\begin{thebibliography}{99}

\bibitem{Co}
CoCoATeam, CoCoA: a system for doing Computations
 in Commutative Algebra, Available at {\tt http://cocoa.dima.unige.it}

\bibitem{ER} J. Eagon, V. Reiner,
Resolutions of Stanley-Reisner rings and Alexander duality.
J. Pure Appl. Algebra  {\bf 130}  (1998) 265--275.

\bibitem{eghp} D. Eisenbud, M. Green, K. Hulek and S. Popescu, Restricting linear
syzygies: algebra and geometry. Compositio Math. {\bf 141} (2005)  1460-1478.

\bibitem{EK} S. Eliahou, M. Kervaire,
Minimal resolutions of some monomial ideals.
J. Algebra {\bf 129} (1990) 1--25.

\bibitem{EV} S. Eliahou, R.H. Villarreal,
The second Betti number of an edge ideal.
XXXI National Congress of the Mexican Mathematical Society
(Hermosillo, 1998), 115--119,
Aportaciones Mat. Comun., {\bf 25},
Soc. Mat. Mexicana, M\'exico, 1999.

\bibitem{faridi:2002} S. Faridi, The facet ideal of a simplicial complex.  Manuscripta Math.
{\bf 109}  (2002)  159--174.

\bibitem{faridi:2004} S. Faridi,
Simplicial trees are sequentially Cohen-Macaulay.
J. Pure Appl. Algebra  {\bf 190}  (2004) 121--136.

\bibitem{Fa} G. Fatabbi, On the resolution of ideals of fat points.
J. Algebra {\bf 242} (2001) 92--108.



\bibitem{FV} C. Francisco, A. Van Tuyl, Sequentially Cohen-Macaulay
edge ideals.  (2005) To appear Proc. Amer. Math. Soc. {\tt math.AC/0511022}

\bibitem{FH} C. Francisco, H. T. H\`a, Whiskers and sequentially Cohen-Macaulay graphs.
(2006) Preprint. {\tt math.AC/0605487}

\bibitem{Fr} R. Fr\"oberg, On Stanley-Reisner rings. In: Topics
in algebra, Banarch Center Publications, {\bf 26} (2) (1990) 57-70.

\bibitem{Gr}  M. Green,
Koszul cohomology and the geometry of projective varieties.
J. Differential Geom. {\bf 19} (1984)  125--171.

\bibitem{HHZ} J. Herzog, T. Hibi, X. Zheng, Cohen-Macaulay chordal graphs.
(2004) Preprint. {\tt math.AC/0407375}

\bibitem{HHZ1} J. Herzog, T. Hibi. X. Zheng, Dirac's theorem on chordal
graphs and Alexander duality. Eur. J. Comb. {\bf 25} (2004) 949-960.


\bibitem{hi} J. Herzog and S. Iyengar, Koszul modules. J. Pure Appl. Algebra
{\bf 201} (2005) 154-188.

\bibitem{Ho} M. Hochster,
Cohen-Macaulay rings, combinatorics, and simplicial complexes.
Ring theory, II (Proc. Second Conf., Univ. Oklahoma, Norman, Okla., 1975),
pp. 171--223. Lecture Notes in Pure and Appl. Math., {\bf 26},
Dekker, New York, 1977.

\bibitem{J} S. Jacques, Betti numbers of graph ideals. Ph.D. Thesis,
University of Sheffield, 2004. {\tt math.AC/0410107}

\bibitem{JK} S. Jacques, M. Katzman, The Betti numbers of forests.
(2005) Preprint. {\tt math.AC/0401226}

\bibitem{K} M. Katzman, Characteristic-independence of Betti numbers
of graph ideals. J. Combinatorial Theory, Series A.
{\bf 113} (2006) 435-454.


\bibitem{MS} E. Miller, B. Sturmfels, {\it Combinatorial Commutative Algebra.}
Springer GTM 227, Springer, 2004.

\bibitem{RVT} M. Roth, A. Van Tuyl, On the linear strand of an edge
ideal. (2004) To appear in Comm. Algebra. {\tt math.AC/0411181}

\bibitem{S} A. Simis,
On the Jacobian module associated to a graph.
Proc. Amer. Math. Soc. {\bf 126} (1998) 989--997.

\bibitem{SVV} A. Simis, W.V. Vasconcelos, R.H. Villarreal,
On the ideal theory of graphs.
J. Algebra {\bf 167} (1994) 389--416.


\bibitem{V} R. H. Villarreal,
Rees algebras of edge ideals.
Comm. Algebra {\bf 23} (1995) 3513--3524.

\bibitem{V1} R. H. Villarreal,
Cohen-Macaulay graphs.
Manuscripta Math. {\bf 66} (1990) 277--293.

\bibitem{V2} R. H. Villarreal,
{\it Monomial algebras.}
Monographs and Textbooks in Pure and Applied Mathematics, {\bf 238}.
Marcel Dekker, Inc., New York, 2001.

\bibitem{Z} X. Zheng, Resolutions of Facet Ideals. Comm. Algebra {\bf 32}
(2004) 2301-2324.

\end{thebibliography}
\end{document}